\documentstyle{article}
\begin{document}
\begin{center}
{\bf SOLUTIONS OF CERTAIN FRACTIONAL KINETIC EQUATIONS 
AND A FRACTIONAL DIFFUSION EQUATION}\\[0.5cm]
{\bf R.K. SAXENA}\\
Department of Mathematics and Statistics, Jai Narain Vyas University\\
Jodhpur-342 004, India\\ [0.5cm]

{\bf A.M. MATHAI}\\
Department of Mathematics and Statistics, McGill University\\
Montreal, Canada H3A 2K6\\
and\\
Centre for Mathematical Sciences, Pala Campus, Pala-686 574, Kerala, India\\[0.5cm]

{\bf H.J.  HAUBOLD}\\
Office for Outer Space Affairs, United Nations\\
P.O. Box 500, A-1400 Vienna, Austria\\
and\\
Centre for Mathematical Sciences, Pala Campus, Pala-686 574, Kerala, India\\
\end{center}
\bigskip
\noindent
{\bf Abstract.} In view of the usefulness and importance of  kinetic equations in certain physical problems, the authors derive the explicit solution of a fractional kinetic equation of general character, that unifies and extends earlier results. Further, an alternative shorter method based on a result developed by the authors is given to derive the solution of a fractional diffusion equation.
 
\section{ Introduction} 

    Fractional reaction/diffusion equations involve fractional derivatives with respect to time and space and  are  studied to describe anomalous reaction/diffusion of dynamic systems with chaotic motion. Fractional kinetic equation for Hamiltonian chaos is discussed by Zaslavsky (1994). Solutions and applications of certain kinetic equations are studied by Saichev and Zaslavsky (1997). Solutions of a fractional kinetic equation is investigated by Haubold and Mathai (2000) for a simple production-destruction mechanism. This equation was generalized by Saxena, Mathai, and Haubold  (2002). In recent articles, Saxena, Mathai, and Haubold (2002, 2004a, 2004b) discussed the solution of a number of generalized fractional kinetic equations . In the present article, we investigate the solution of a unified fractional kinetic equation, which provides unification and extension of  results on fractional kinetic equations given earlier by Haubold and Mathai (2000) and Saxena, Mathai, and Haubold  (2002, 2004a). We also present the solution of a fractional integral equation discussed by Miller and Ross (1993). Further, an alternative proof of the solution of a fractional diffusion equation given earlier by Kochubei (1990) is investigated, which is based upon a result given by Saxena, Mathai, and Haubold  (2006). Most of the results are obtained in terms of generalized Mittag-Leffler functions in  elegant and compact forms, which are suitable for numerical computation.
    
     The paper is organized as follows. Section 2 contains the solution of a unified fractional kinetic equation while Section 3 considers special cases of the equation. A shorter alternative method for the solution of a diffusion equation discussed earlier by Kochubei (1990) is presented in Section 4. A series representation and asymptotic expansion of the solution are given in Section 5. Incidentally, an H-function representation of a one-sided L\'{e}vy stable density is also obtained.

\section{ Unified fractional kinetic equations}

    In this Section, we present a method based on Laplace transform for deriving the solution of the unified fractional kinetic equations. \par
\medskip
\noindent
{\bf Theorem 1.} If $Re(\nu_j)>0, a_j>0, j\in \mbox{N}$, and  f(t) be a given function,  defined on $\Re_+$, then the equation
\begin{equation}
N(t)-N_0f(t)=-\sum^n_{j=1}a_j\;_0D_t^{-\nu_j} N(t),
\end{equation}
is solvable  and its particular solution is given by
\begin{eqnarray}
N(t)&=& N_0\sum^\infty_{l=0}(-1)^l \sum_{r_1+\ldots+r_{n-1}=l}\frac{(l)!}{(r_1)!\ldots(r_{n-1})!}\left\{\prod^{n-1}_{\mu=1}(a_{\mu+1})^{r_\mu}\right\}\nonumber\\
&&\int_0^t f(u)(t-u)^{\sum^{n-1}_{\mu=1}\nu_{\mu+1}-1}E^{(l+1)}_{\nu_1,\sum^{n-1}_{\mu=1}\nu_{\mu+1}}[-a_1(t-u)^{\nu_1}]du,
\end{eqnarray}
where the summation in (2) is taken over all nonnegative integers $r_1,\ldots,r_n$  such that $r_1+\ldots+r_{n-1}=l$, and provided that the series and integral in (2)  are convergent. Here $_0D_t^{-\nu_j}, j\in \mbox{N}$  are Riemann-Liouville fractional integrals, defined by 
\begin{equation}
_0D_t^{-\nu} f(t)=\frac{1}{\Gamma(\nu)}\int_0^t (t-u)^{\nu-1}f(u)du, Re(\nu) >0,
\end{equation}
with $_0D_t^0 f(t)=f(t)$  (Oldham and Spanier, 1974;  Miller and Ross, 1993), $E^\delta_{\beta, \gamma}(z)$  is the generalized Mittag-Leffler function,  defined by  Prabhakar (1971) in terms of series representation  as 
\begin{equation}
E^\delta_{\beta, \gamma}(z)=\sum^\infty_{\tau=0}\frac{(\delta)_\tau z^\tau}{\Gamma(\beta\tau+\gamma)(\tau)!}\;\;(\beta, \gamma, \delta\in C, Re(\beta)>0, Re(\gamma)>0).
\end{equation}
\noindent
{\bf Proof.} By the application of the convolution theorem of the Laplace transform (Erd\'{e}'lyi et al., 1953, p. 259)  to  (3), we find that
\begin{eqnarray}
L\left\{_0D_t^{-\nu}f(t);s\right\} &=&L\left\{\frac{t^{\nu-1}}{\Gamma(\nu)}\right\}L( f(t)),\nonumber\\
&=& s^{-\nu}f^\sim(s),
\end{eqnarray}
where $f^\sim(s)=\int_0^\infty e^{-st} f(t)dt, s\in C, Re(s)>0$. 
Applying Laplace transform to (1) and using (5), it gives
\begin{eqnarray}
N^\sim(s)&=&\frac{N_0 f^\sim(s)}{1+a_1 s^{-\nu_1}+\ldots+ a_n s^{-\nu_n}}\\
&=& N_0 f^\sim (s)\sum^\infty_{l=0}(-1)^l\frac{\left(\sum^{n-1}_{j=1}a_{j+1}s^{-\nu_{j+1}}\right)^l}{(1+a_1 s^{-\nu_1})^{l+1}}, \left|\frac{\sum^n_{j=2}a_j s^{-\nu_j}}{1+a_1 s^{-\nu_1}}\right|<1.\nonumber
\end{eqnarray}
  
If we employ the identity (Abramowitz and Stegun, 1968, p. 823)
\begin{equation}
(x_1+\ldots +x_m)^l=\sum_{r_1+\ldots +r_n=l}\frac{(l)!}{(r_1)!\ldots (r_n)!}\prod^m_{\mu=1} x_\mu ^{r_\mu},
\end{equation}	
where the summation is taken over all nonnegative integers, $r_1, \ldots, r_n$, such that $r_1+\ldots +r_n=l$, then for $|a_1s^{-\nu_1}|<1$, (7) transforms into the form
\begin{eqnarray}
N^\sim(s)&=& N_0 f^\sim(s)\sum^\infty_{l=0}(-1)^l\sum_{r_1+\ldots+r_{n-1}=l\atop
r_1>\ldots r_{n-1}>0}\frac{(l)!}{(r_1)!\ldots(r_{n-1})!}\nonumber\\
&&\frac{\left\{\prod^{n-1}_{\mu=1}(a_{\mu+1})^{r_\mu}\right\}s^{-\sum^{n-1}_{\mu=1}\nu_{\mu+1}}}{(1+a_1 s^{-\nu_1})^{l+1}}.
\end{eqnarray} 
Taking the inverse Laplace transform of (8) by making use of  the formula (Kilbas, Saigo, and Saxena, 2004, eq. (12)) 
\begin{equation}
L^{-1}\left\{s^{-\gamma}(1-as^{-\beta})^{-\delta}; t\right\}= t^{\gamma-1}E^\delta_{\beta,\gamma}(at^\beta),
\end{equation}         
where 
$Re(s)>|a|^{1/Re(\gamma)}, Re(\gamma)>0, Re(s)>0,$\\
and applying the convolution theorem of the Laplace transform, the result (2) is established.\par
\medskip
\noindent
{\bf Remark 1.} The generalized Mittag-Leffler function defined by (4) is studied by Prabhakar (1971) and  Kilbas, Saigo, and Saxena (2004). Recently this function is used in the theory of finite-size scaling of systems with strong anisotropy and long-range interaction by Chamati and Tonshev  (2006).\par
\medskip
\noindent
\section  {Special cases}

    Some special cases of  Theorem 1 are of interest to be highlighted.
If we set $\nu_j=j\nu, a_j=(^n_j)c^{j\nu}(j\in \mbox{N})$, we obtain\par
\medskip
\noindent
{\bf Theorem 2.} If $Re(\nu)>0, c>0$  and $ f(x)\in \Re_+$, then the equation
\begin{equation}
N(t)-N_0f(t)=-\sum^n_{r=1}(^n_r) c^{\nu r} D_t^{-\nu r} N(t),
\end{equation}
is solvable and its solution has the form
\begin{equation}
N(t)=N_0\frac{d}{dt}\int_0^t f(u)E^n_{\nu,1}[-c^\nu(t-u)^\nu]du,
\end{equation}  
where $E^n_{\nu,1}(x)$ is the generalized Mittag-Leffler function defined by (4) and provided that the integral (11) is convergent.

     When  $n=1$, we obtain the following result given by Hille and Tamarkin (1930).\par
\medskip
\noindent
{\bf Corollary 2.1.} Let $Re(\nu)>0, c>0$  and let $f(x)\in \Re_+$, then  for the solution of the integral equation
\begin{equation}
N(t)-N_0f(t)=-c^\nu\;_0D_t^{-\nu} N(t),
\end{equation}
holds the following formula 
\begin{equation}
N(t)=N_0\frac{d}{dt}\int_0^tf(u)E_\nu[-c^\nu (t-u)^\nu ]du,
\end{equation}
where $E_\nu(z)$  is an entire function of  order $\rho=\frac{1}{\nu}$   and type $\sigma=1$, defined by 
\begin{equation} 
E_\nu(z)=\sum^\infty_{\mu=1}\frac{z^\mu}{\Gamma(\mu \nu+1)}, \;(\nu\in C, Re(\nu)>0).
\end{equation}
\noindent 
{\bf Note 1.} The above result has also been given by the authors in a different form (Saxena, Mathai, and Haubold,  2004a, 2004b).
 
     If we set $f(t)=t^{\gamma-1}E^\delta_{\nu,\gamma}[-(ct)^\nu]$,  Theorem 2 yields\par
\medskip
\noindent
{\bf Corollary 2.2.}  Let $Re(\nu)>0, Re(\gamma)>0, c>0$, then for the solution of the kinetic equation 
\begin{equation}
N(t)-N_0t^{\gamma-1}E_{\nu,\gamma}^\delta[-(ct)^\nu]=-\sum^n_{r=1}(^n_r) c^{r\nu}\;_0D_t^{-r\nu}n\in \mbox{N}
\end{equation}
holds the relation
\begin{equation} 
N(t)=N_0t^{\gamma-1}E^{\delta+n}_{\nu, \gamma}[-(ct)^\nu], n\in \mbox{N}.
\end{equation}     
For $f(t)=t^{\rho-1}$, Theorem 2 yields the following result \par
\medskip
\noindent
{\bf Corollary 2.3.} If $Re(\rho)>0, Re(\nu)>0, c>0$, then for the solution of the equation
\begin{equation}
N(t)-N_0t^{\rho-1}=-\sum^n_{r=1}(^n_r)c^{r\nu}\;_0D_t^{-r\nu}N(t), r\in \mbox{N},
\end{equation}
holds the relation
\begin{equation}
N(t)=N_0t^{\rho-1}\Gamma(\rho)E^n_{\nu, \rho}[-(ct)^\nu], r\in \mbox{N}.
\end{equation}
For $n=1$, eq. (18) reduces to a result given by Saxena, Mathai, and Haubold  (2002, p. 283, eq. (15)). When $a_j=a^js^{\nu j}$, for $j=1,\ldots, n$, we obtain\par
\medskip
\noindent
{\bf Theorem 3.} Let $Re(\nu)>0, a>0, t>0, n>1, |a^{n+1}s^{-(n+1)\nu}|<1$, and $f(x)$ be a given function defined on $\Re_+$, then the equation
 \begin{equation}
N(t)-N_0f(t)=-\sum^n_{r=1}a^r\;_0D_t^{-\nu r} N(t),
\end{equation}
is solvable and its solution is given by
\begin{eqnarray}
N(t)&=& N_0\left\{\frac{d}{dt}\int_0^t f(u)E_{(n+1)\nu,\nu}[a^{n+1}(t-u)^{(n+1)\nu}]du\right.\nonumber\\
&&\left.-a \int_0^t(t-u)^{\nu-1}E_{(n+1)\nu,\nu}[a^{n+1}(t-u)^{(n+1)\nu}]du\right\},
\end{eqnarray}
where $E_{(n+1)\nu,\nu}(z)$    is the generalized Mittag-Leffler function $E_{\alpha,\beta}(z)$ defined as
\begin{equation}
E_{\alpha, \beta}(z)=\sum^\infty_{\mu=1}\frac{z^\mu}{\Gamma(\alpha \mu+\beta)},\;(\alpha, \beta\in C, Re(\alpha)>0, Re(\beta)>0)
\end{equation}
and provided that the integral in  (20) is convergent.                                

      If we take $\nu^j=j\nu$,   for $j=1,\ldots n$,  then it is interesting to note that Theorem 1 yields the following result given by (Miller and Ross, 1993) in a different form:\par
\medskip
\noindent
{\bf Theorem 4.}  Let $Re(\nu)>0, a_j>0,$ and $f(x)$ be a given function defined on $\Re_+, |a_1 s^{-\nu}|<1$, then the fractional  kinetic equation 
\begin{equation}
N(t)-N_0f(t)=-\sum^n_{j=1}a_j\;_0D_t^{-j\nu} N(t),
\end{equation}
is solvable and has the  solution given by 
\begin{eqnarray}
N(t)&=& N_0\sum^\infty_{l=0}(-1)^l \sum_{r_1+\ldots+r_{n-1}=l}\frac{(l)!}{(r_1)!\ldots(r_{n-1})!}\left\{\prod^{n-1}_{\mu=1}(a_{\mu+1})^{r_\mu}\right\}\\
&\times& \int_0^t f(u)(t-u)^{\sum^{n-1}_{\mu=1}\nu(\mu+1)r_\mu-1}E^{(l+1)}_{\nu_1,\sum^{n-1}_{\mu=1}\nu(\mu+1)r_\mu}[-a_1(t-u)^{\nu_1}]du,\nonumber
\end{eqnarray}
provided that the series and integral in (23) are convergent.
\section{Fractional diffusion equation}

    In this Section we present an alternative shorter method for deriving the solution of a diffusion equation discussed earlier by Kochubei (1990).\par
\medskip
\noindent
{\bf Theorem 5.} Consider the Cauchy problem
\begin{equation}
_0D_t^\alpha N(x,t)=-c^\nu\Delta N(x,t),\;(0<\alpha<1; x\in \Re^n; \;0<t\leq T),
\end{equation}
with
\begin{equation}
N(x,t=0)=\delta(x), x\in \Re, \lim_{x\rightarrow \pm \infty} N(x,t)=0
\end{equation}
$_0D_t^\alpha$  is the regularized  Caputo (1969)  partial fractional derivative with respect to t, defined by
$$_0D_t^{\alpha}N(x,t)=\frac{1}{\Gamma(1-\alpha)}\left[\frac{\partial}{\partial t}\int_0^t\frac{N(x,s)ds}{(t-x)^\alpha}-\frac{N(x,0)}{t^\alpha}\right],$$
and $\Delta$  is the Laplacian. The fundamental solution of the above Cauchy problem is given by 
\begin{equation}
N(x,t)=|x|^{-n}\pi^{-\frac{n}{2}}H^{2,0}_{1,2}\left[\frac{|x|^2t^{-\alpha}}{4c^\nu}|^{(1,\alpha)}_{(n/2,1),(1,1)}\right],
\end{equation}
where $H^{2,0}_{1,2}(.)$    is the H-function (Mathai and Saxena, 1978).\par
\medskip
\noindent
{\bf Proof.} Applying the Laplace transform with respect to t, using the result (Caputo, 1969)
$$L\left\{_0D_t^\alpha  N(x,t)\right\}=s^\alpha N^\sim(x,s)-\sum^{m-1}_{r=0}s^{\alpha-r-1} N^{(r)}(x,0),$$
$$m-1<\alpha\leq m,\;m\in \mbox{N},$$
and Fourier transform with respect to $x$, gives
$$s^\alpha N^{\sim^*}(k,s)-s^{\alpha-1}=-c^\nu|k|^2 N^{\sim^ *} (k,s),$$
where the symbol ''$\sim$`` indicates the Laplace transform with respect to  the time variable $t$ and the symbol ''*`` the Fourier transform with respect to the space variable $x$.\par
\noindent
Solving for $N^{\sim^*}$, we have
\begin{equation}
N^{\sim^*}(k,s)=\frac{s^{\alpha-1}}{s^\alpha+c^\nu|k|^2}.
\end{equation}
By virtue of the following Fourier  transform formula (Samko, Kilbas, and Marichev, 1990, p. 538, eq. (27.1))

\begin{equation}
\left(F_x\left[|x|^{(2-n)/2}K_{(n-2)/2}(a|x|)\right]\right)(\tau)=\left(\frac{2\pi}{a}\right)^{n/2}\frac{a}{a^2+\tau^2}, 
(\tau\in \Re^n; n\in \mbox{N}, a>0),
\end{equation}
where the multidimensional Fourier transform with respect to $x\in \Re^n$ is defined by 
\begin{equation}
(F_x N)(\tau, t)=\int_{\Re^n}N(x,t)e^{ix \tau} dx\; (\tau\in \Re^n; t>0)
\end{equation}             
and $K_\nu(.)$ is the  modified Bessel function of the second kind,  yields
\begin{equation}
\tilde{N}(x,s)=c^{-\nu}s^{\alpha-1}(2\pi)^{-\frac{n}{2}}\left(\frac{|x|c^{\nu/2}}{s^{\alpha/2}}\right)^{1-\frac{n}{2}}K_{n-2/2}\left[\frac{|s^{\frac{\alpha}{2}}|x|}{c^{\frac{\nu}{2}}}\right].
\end{equation}

     In order to invert the Laplace transform, we employ the following result given by the authors (Saxena, Mathai, and Haubold, 2006)
\begin{equation}
L^{-1}\left\{s^{-\rho}K_\nu(zs^\sigma);t\right\}=\frac{1}{2}t^{\rho-1}H^{2,0}_{1,2}\left[\frac{z^2t^{-2\sigma}}{4}\left|
^{(\rho, 2\sigma)}_{(\frac{\nu}{2},1)(-\frac{\nu}{2},1)}\right.\right],
\end{equation}
where $K_\nu(x)$ is the modified Bessel function of the second kind, $Re(z^2)>0, Re(s)>0$. Thus we obtain the solution in the form       
\begin{eqnarray}
N(x,t)&=&\frac{1}{2}(2\pi)^{-\frac{n}{2}}c^{-\frac{\nu}{2}-\frac{n\nu}{4}}|x|^{1-\frac{n}{2}}t^{-\frac{\alpha}{2}-\frac{\alpha n}{4}}\nonumber\\
&& H^{2,0}_{1,2}\left[\frac{t^{-\alpha}|x|^2}{4c^\nu}\left|^{(1-\frac{\alpha}{2}-\frac{\alpha n}{4}, \alpha)}_{(\frac{n-2}{2},1),(\frac{2-n}{4},1)}\right]\right..
\end{eqnarray}
By virtue of a result in Mathai and Saxena (1978),
\begin{equation}
x^\sigma H^{m,n}_{p,q}\left[x|^{(a_p, a_p)}_{(b_q, b_q)}\right]=H^{m,n}_{p,q}\left[x|^{(a_p+\sigma A_p, A_p)}_{(b_q+\sigma B_q, B_q)}\right],
\end{equation}
the power of  the expression $[\left\{t^{-\nu}|x|^2\right\}/4c^\nu]$ can be absorbed inside the H-function and consequently we obtain
\begin{equation}
     N(x,t)  = |\pi^{\frac{1}{2}}x|^{-n}\;H^{2,0}_{1,2}\left[\frac{t^{-\alpha}|x|^2}{4c^\nu}\left|^{(1,\alpha)}_{(\frac{n}{2},1),(1,1)}\right]\right..
\end{equation}   
\noindent
{\bf Remark 1.} If we employ the identity  (Mathai and Saxena, 1978)
\begin{equation}
H^{m,n}_{p,q}\left[ x^\lambda|^{(a_p, A_p)}_{(b_q, B_q)}\right]=\frac{1}{\lambda}H^{m,n}_{p,q}\left[x|^{(a_p,A_p/\lambda)}_{(b_q, B_q/\lambda)}\right],\lambda>0
\end{equation}   
the solution given by  (32) can be expressed in the form 
\begin{equation}
 N(x,t)  =  \frac{1}{\alpha}|\pi^{\frac{1}{2}}x|^{-n}\;H^{2,0}_{1,2}\left(\frac{|x|^2}{4c^\nu t^\alpha}\right)^{1/\alpha}\left| ^{(1,1)}_{(\frac{n}{2}, \frac{1}{\alpha}), (1,\frac{1}{\alpha})}\right.,
\end{equation}
where $\alpha>0$. We also note that the above form of the solution is due to Schneider and Wyss (1989). There is one importance of our result (32) that it includes the L\'{e}vy stable density in terms of the H-function as shown in (45).  Similarly, using  the identity  (35) we arrive at  
\begin{equation}
N(c,t)=\frac{1}{2}|\pi^{\frac{1}{2}}x|^{-n}H^{2,0}_{1,2}\left[\frac{t^{-\frac{\alpha}{2}}|x|}{2c^{\frac{\nu}{2}}}\left|^{(1,\frac{\alpha}{2})}_{(\frac{n}{2},\frac{1}{2}), (1,\frac{1}{2})}\right]\right.,
\end{equation}
where $n$ is not an even integer.
     This form of the H-function is useful in determining its expansion in powers of $x$. Due to importance of the solution, we also discuss its series representation and behavior.  
 
\section{Series representation of the solution}
 
    Using the series expansion for the H-function given in Mathai and Saxena (1978), it follows  that 
\begin{eqnarray}
&&H^{2,0}_{1,2}\left[x\left|^{(1,1)}_{(\frac{n}{2}, \frac{1}{\alpha}),(1,\frac{1}{\alpha})}\right]=\frac{1}{2\pi i}\int_L\frac{\Gamma(\frac{n}{2}-\frac{s}{\alpha})\Gamma(1-\frac{s}{\alpha})}{\Gamma(1-s)}x^s ds\right.\\
&=&\alpha\left\{\sum^\infty_{l=0}\frac{\Gamma(1-\frac{n}{2}-l)(-1)^lx^{\alpha(\frac{n}{2}+l)}}{\Gamma(1-\frac{a n}{2}-\alpha l)(l)!}+\sum^\infty_{l=0}\frac{\Gamma(\frac{n}{2}-1-l)(-1)^lx^{\alpha(1+l)}}{\Gamma(1-\alpha-\alpha l)(l!)}\right\},\nonumber
\end{eqnarray}    
where $n$ is not an even integer.

     Thus for $n=1$, we find that
 \begin{equation}    
N(x,t) = \frac{1}{2t^{\frac{\alpha}{2}}}\sum^\infty_{l=0}(-1)^l\frac{A^{\frac{l}{2}}}{\Gamma(1-\alpha(l+1)/2)(l!)}, 
\end{equation}
where  $A=\frac{x^2}{t^\alpha}$  and the duplication formula for the gamma function is used. 

For $n=2$, the H-function of (37)  is singular and in this case, the result  is explicitly given by  Barkai (2001) in the form 
\begin{equation}
N(x,t)\sim \frac{1}{\pi\Gamma(1-\alpha)t^\alpha}ln[\frac{t^{\frac{\alpha}{2}}}{x}].
\end{equation}
For $n=3$, the series expansion is given by
\begin{equation}
N(x,t)=\frac{1}{4\pi t^{3\alpha/2}A^{1/2}}\sum^\infty_{l=0}\frac{(-1)^l A^{l/2}}{\Gamma[1-\alpha(1+l/2)]}.
\end{equation}

     From above, it readily follows that for $n=3$ and $\alpha \neq 1$
\begin{equation}
N(x,t)\sim\frac{1}{x},\; \mbox{as}\;\; x\rightarrow \infty.
\end{equation}
    
        It will not be out of place to mention that the one sided L\'{e}'vy stable density can be obtained from Laplace inversion formula (31) by virtue of the identity
\begin{equation}
K_{\pm\frac{1}{2}}(x)=\left(\frac{\pi}{2x}\right)^{\frac{1}{2}} e^{-x},
\end{equation}
and can be conveniently expressed in terms of the Laplace transform 

      \begin{equation}
\int_0^\infty e^{-ut}\Phi_\rho(t)dt=e^{-u^\rho},\;\;Re(u)>0,\; Re(\rho)>0.
\end{equation}
The result is
\begin{equation}
\Phi_\rho(t)=\frac{1}{\rho}H^{1,0}_{1,1}\left[\frac{1}{t}\left|^{(1,1)}_{(\frac{1}{\rho}, \frac{1}{\rho})}\right.\right], \;(\rho>0).
\end{equation}
{\bf Note 2.} This result is obtained earlier by Schneider and Wyss (1989) by following a different procedure. Asymptotic behavior of $\Phi_\alpha(t)$   is also given by Schneider (1986). 

      In conclusion, we mention that some of the results derived in this article may find some applications in problems associated with models of long-memory processes driven by L\'{e}'vy noise and other related problems, see the article by Anh, Heyde, and Leonenko  (2002).\par
\bigskip
\noindent
{\bf Acknowledgment.} The authors would like to thank the Department of Science and Technology, Government of India, New Delhi, for the financial assistance for this work under project No. SR/S4/MS:287/05 which enabled this collaboration possible.\par
\bigskip
\noindent
{\bf References}\par
\medskip
\noindent
Anh, V.V., Heyde, C.C., and Leonenko, N.N.:  2002, Dynamic models driven by L\'{e}vy 
noise, {\it Journal of Applied  Probability}, {\bf 39}, 730-747.\par
\smallskip
\noindent  
Abramowitz, M. and Stegun, I.A.: 1968, {\it Handbook of Mathematical Functions with
Formulas, Graphs, and Mathematical Tables}, Dover Publications, Inc. New York.\par
\smallskip
\noindent  
Barkai, E.: 2001, Fractional Fokker-Planck equation, solution, and application,
{\it Physical Review} E, {\bf 63}, 046118.\par
\smallskip
\noindent  
Caputo, M.: 1969, {\it Elasticita  e Dissipazione}, Zanichelli, Bologna.\par
\smallskip
\noindent  
Chamati, H. and Tonchev, N.S.: 2006, Generalized Mittag-Leffler functions in the  theory of  finite-size scaling for systems with strong anisotropy and/or long-range interaction, {\it Journal of Physics A. Mathematical and General}, {\bf 39}, 469- 470.\par
\smallskip
\noindent
Erd\'{e}lyi, A., Magnus, W., Oberhettinger, F., and Tricomi, F.G.: 1953, {\it Higher Transcendental Functions} , Vol. {\bf 1}, McGraw-Hill, New York, Toronto, and London.\par
\smallskip
\noindent 
Haubold, H.J. and Mathai, A.M.: 2000, The fractional kinetic equation and thermonuclear functions, {\it Astrophysics and Space Science}, {\bf 273}, 53-63.\par
\smallskip
\noindent     
Hille, E. and Tamarkin, J.D.: 1930, On the theory of linear integral equations, {\it Annals of Mathematics}, {\bf 31}, 479-528.\par
\smallskip
\noindent        
Kochubei, A.N.: 1990, Diffusion of fractional order, {\it Differential Equations}, {\bf 26}, 485-492.\par
\smallskip
\noindent           
Kilbas, A.A., Saigo, M. and Saxena, R.K.: 2004, Generalized Mittag-Leffler function and generalized fractional calculus operators, {\it Integral Transforms and Special Functions}, {\bf 15}, 31-49.\par
\smallskip
\noindent           
Mathai, A.M. and Saxena, R.K.: 1978, {\it The H-function with Applications in Statistics and Other Disciplines}, John Wiley and Sons, Inc., New York, London, and Sydney. \par
\smallskip
\noindent              
Miller, K.S. and Ross, B.: 1993, {\it An Introduction to the Fractional Calculus and Fractional Differential Equations},  John Wiley and Sons, New York.\par
\smallskip
\noindent               
Oldham, K.B. and Spanier, J.: 1974, {\it The Fractional Calculus. Theory and Applications of Differentiation and Integration of Arbitrary Order}, Academic Press, New York \par
\smallskip
\noindent              
Prabhakar, T.R.: 1971, A singular integral equation with a generalized Mittag-Leffler function in the kernel, {\it Yokohama Journal of Mathematics}, {\bf 19}, 7-15.\par
\smallskip
\noindent              
Saichev, A.I. and Zaslavsky, G.M.: 1997, Fractional kinetic equations: solutions and applications, {\it Chaos}, {\bf 7}, 753-784.\par
\smallskip
\noindent           
Samko, S.G., Kilbas, A.A. and Marichev, O.I.: 1990, {\it Fractional Integrals and Derivatives: Theory and Applications}, Gordon and Breach Science Publishers, New York.\par
\smallskip
\noindent           
Saxena, R.K., Mathai, A.M., and Haubold, H.J.: 2002, On fractional kinetic equations, {\it Astrophysics and Space Science}, {\bf 282}, 281-287. \par
\smallskip
\noindent           
Saxena, R.K., Mathai, A.M., and Haubold, H.J.: 2004a, On generalized fractional kinetic equations, {\it Physica  A},  {\bf 344}, 657-664.\par
\smallskip
\noindent           
Saxena, R.K., Mathai, A.M. and Haubold, H.J.: 2004b, Unified fractional kinetic equation and a fractional diffusion equation, {\it Astrophysics and Space Science}, {\bf 290}, 299-310.\par
\smallskip
\noindent           
Saxena, R.K., Mathai, A.M., and Haubold, H.J.: 2006, Solution of generalized fractional reaction-diffusion equations, {\it Astrophysics and Space Science}, {\bf 305}, 305-313.\par
\smallskip
\noindent           
Schneider, W.R.: 1986, in {\it Stochastic Processes in Classical and Quantum Systems}, S. Albeverio, G. Casati, and D. Merlini (Eds.), Springer-Verlag, Berlin.\par
\smallskip
\noindent           
Schneider, W.R. and Wyss, W.: 1989,  Fractional diffusion and wave equation, {\it Journal of Mathematical Physics}, {\bf 30}, 134-144. \par
\smallskip
\noindent                 
Zaslavsky, G.M.: 1994, Fractional kinetic equation for Hamiltonian chaos, {\it Physica D}, {\bf 76}, 110-122.
\end{document}